\documentclass[12pt]{article}          

\usepackage{theorem}
\usepackage{oldgerm}
\usepackage{amsmath}
\usepackage{amssymb}
\usepackage{eclbip}
\usepackage{amscd}

{\theorembodyfont{\rmfamily}

                             }

\theoremstyle{break} 
\def\preprint{}
\def\finished{\begin{tabbing}
Keywords: \={\tt 14J32} Calabi-Yau manifolds, mirror symmetry\\
\>{\tt14M25} Toric varieties, Newton polyhedra\\
\>{\tt81T30} String and superstring theories
\end{tabbing}}
\def\archive {math.AG/0103214}

\def\Title{     Toric complete intersections and \\[3mm]
                        weighted projective space
}

\long\def\Abstract{
It has been shown by Batyrev and Borisov that nef partitions of reflexive
polyhedra can be used to construct mirror pairs of complete intersection
Calabi--Yau manifolds in toric ambient spaces. We construct a number of such
spaces and compute their cohomological data. We also discuss the relation of
our results to complete intersections in weighted projective spaces and try
to recover them as special cases of the toric construction.
As compared to hypersurfaces, codimension two more than doubles the number
of spectra with $h^{11}=1$. Altogether we find 87 new
(mirror pairs of) Hodge data, mainly with $h^{11}\le4$.
}

\usepackage{cite,isolatin1,amssymb,latexsym}    \def\bye{\end{document}}
\long\def\new#1\endnew{{\bf #1}}                \long\def\del#1\enddel{}
\def\ifundefined#1{\expandafter\ifx\csname#1\endcsname\relax}

\def\BC{\begin{center}} \def\EC{\end{center}}   \catcode`\"=\active \let"=\"

\def  \b{\begin}
\def  \e{\end}

\def \d{\hat}
\def \beq{\b{equation}}
\def \eeq{\e{equation}}
\def \bea{\b{eqnarray}}
\def \eea{\e{eqnarray}}

\def \pro{{\em Proof:\ }}

\def \R{\Bbb{R}}
\def \Z{\Bbb{Z}}
\def \P{\Bbb{P}}

\let\0=\over      \let\1=\vec      \def\2{{1\over2}}
\let\4=\underline \let\5=\overline \let\6=\partial      \def\7#1{{#1}\llap{/}}
\def\8#1{{\textstyle{#1}}}         \def\9#1{{\ifmmode{\pmb{#1}}\else\bf#1\fi}}

\def\eeql#1 {\label{#1}\eeq}      
\def\beq{\begin{equation}}      \def\eeq{\end{equation}}
\def\bea{\begin{eqnarray}}      \def\eea{\end{eqnarray}}

\def\EEL#1 {\label{#1}\EE}           
\def\BE {\begin{equation}}      \def\EE {\end{equation}}
\def\BEA{\begin{eqnarray}}      \def\EEA{\end{eqnarray}}

\def\mao#1{\mathop{\rm #1}\nolimits}

\let\and=\wedge

\let\bra=\langle        \let\ket=\rangle        \def\<#1\>{\bra #1 \ket}
       \def\rel#1 #2{\buildrel #1 \over {#2}}

\def\fnote#1#2{\begingroup\def\thefootnote{#1}\footnote{#2}
                \addtocounter{footnote}{-1}\endgroup}

    \let\h=\eta     
   \let\l=\lambda

               \let\S=\Sigma
       \let\D=\Delta

\def\IR{{\mathbb R}} \def\IC{{\mathbb C}} \def\IP{{\mathbb P}}
  
\def\IZ{{\mathbb Z}}

\def\plb#1 #2 {Phys. Lett. {\bf B#1} #2 }
\def\phr#1 #2 {Phys. Rep. {\bf  #1} #2 }        
\def\npb#1 #2 {Nucl. Phys. {\bf B#1} #2 }
\def\aph#1 #2 {Ann. Phys. {\bf #1} #2 }         
\def\jmp#1 #2 {J. Math. Phys. {\bf #1} #2 }
\def\jgp#1 #2 {J. Geom. Phys. {\bf #1} #2 }
\def\prd#1 #2 {Phys. Rev. {\bf D#1} #2 }
\def\prl#1 #2 {Phys. Rev. Lett. {\bf #1} #2 }
\def\rmp#1 #2 {Rev. Mod. Phys.  {\bf #1} #2 }
\def\zpc#1 {Z. Phys. {\bf #1C} }
\def\cmp#1 #2 {Commun. Math. Phys. {\bf #1} #2 }
\def\cqg#1 #2 {Class.Quant.Grav. {\bf #1} #2 }
\def\mpl#1 {Mod. Phys. Lett. {\bf A#1} }
\def\cpc#1 {Computer Phys. Commun. {\bf #1} }   
\def\ijmp#1 {Int. J. Mod. Phys. {\bf A#1} }
\def\ijmpC#1 {Int. J. Mod. Phys. {\bf C#1} }
\def\atmp#1 #2 {Adv. Theor. Math. Phys. {\bf #1} #2 }

\def\BP{\begin{picture}} \def\EP{\end{picture}}         
\newcounter{TRefNX} \let\OLDcite=\cite  \makeatletter
\def\makeTRefs#1{\@for  \NewTRef:=#1\do{\global\makeTRef{\NewTRef}}}
\def\makeTRef#1{\ifundefined{TRef#1}\stepcounter{TRefNX}%
\expandafter\xdef\csname TRef#1\endcsname{\theTRefNX}\fi}\makeatother
\def\NEWcite#1{\makeTRefs{#1}\OLDcite{#1}}

\ifundefined{draftmode} {} \else        \let\cite=\NEWcite
   \newcount\HOUR  \HOUR=\the\time \divide\HOUR by 60    \multiply \HOUR by 60
   \newcount\MIN   \MIN=\the\time  \advance\MIN by -\HOUR   \divide\HOUR by 60
   \ifnum   \MIN>9  \def\printTIME{{\it\the\HOUR\,:\,\the\MIN}}
         \else   \def\printTIME{{\it\the\HOUR\,:\,0\the\MIN}} \fi 
   
   \def\LLab#1{\BP(0,0)\unitlength=1mm\put(-12,.5){\makebox(0,0)[cr]{\small #1
        \rlap{$_{_{\makeatletter\csname TRef#1\endcsname\makeatother}}$}}}\EP}
\fi

\begin{document}

\newcommand{\Int}{\operatorname{Int}}
\newcommand{\Co}{\operatorname{Conv}}
\newcommand{\m}{\operatorname{Specm}}
\newcommand{\Div}{\operatorname{Div}}
\newcommand{\Pic}{\operatorname{Pic}}
\newcommand{\SF}{\operatorname{SF}}
\newcommand{\ELB}{\operatorname{ELB}}
\newcommand{\DivT}{\operatorname{Div_T}}
\newcommand{\PDiv}{\operatorname{PDiv}}
\newcommand{\pdiv}{\operatorname{div}}
\newcommand{\CDiv}{\operatorname{CDiv}}
\newcommand{\CDivT}{\operatorname{CDiv_T}}
\newcommand{\Hom}{\operatorname{Hom}}
\newcommand{\ke}{\operatorname{kern}}
\newcommand{\di}{\operatorname{dim}}
\newcommand{\Proj}{\operatorname{Proj}}

\vspace*{-18pt}\begin{flushright}  \archive\\[3pt] \preprint  \end{flushright}
\vspace{-3mm}

\BC{\huge\bf \Title}
\\[9mm]
        Maximilian KREUZER,\fnote{\#}{e-mail: kreuzer@hep.itp.tuwien.ac.at}
        Erwin RIEGLER\fnote{$\,\Box$\,}{e-mail: riegler@hep.itp.tuwien.ac.at}
\\[5mm]
        Institut f\"ur Theoretische Physik, Technische Universit\"at Wien\\
        Wiedner Hauptstra\ss e 8--10, A-1040 Wien, AUSTRIA
\\[5mm]                       and
\\[5mm] David A. SAHAKYAN\fnote{\,*\,}{e-mail: sahakian@theory.uchicago.edu}
\\[3mm] Department of Physics, University of Chicago\\
        5640 S. Ellis Av., Chicago, IL 60637, USA
\\[7mm]
\vfil
{\bf ABSTRACT } \\[7mm]  \parbox{15cm}
{\baselineskip=14.5pt ~~~\Abstract} \EC
\vfil \noindent \preprint \\ \finished 
\setcounter{page}{0} \thispagestyle{empty} \newpage \pagestyle{plain}

\section{Introduction}

The first sizeable sets of Calabi--Yau manifolds were constructed as complete
intersections (CICY) in products of projective spaces \cite{gr87,CICY}. These
manifolds have many complex structure deformations
but only few Kähler moduli, which are inherited from the ambient space.
With the discovery of mirror symmetry \cite{lvw} the main interest therefore
turned to weighted projective ($W\P$) spaces, where the resolution of
singularities contributes additional Kähler moduli and thus provides
a much more symmetric picture \cite{CLS}. It turned out, however, that
mirror symmetry is only approximately realized in this class of models
\cite{nms,KlS}.

It was then discovered by Batyrev \cite{bat} that toric geometry (TG), which
generalizes (products of) $W\P$ spaces, provides the appropriate framework
for mirror symmetry: In TG the monomial deformations of the hypersurface
equations and the gluing data defining the ambient space are given in
terms of lattice polytopes that live in a dual pair of lattices.
The Calabi--Yau condition for the generic hypersurface requires that these
polytopes are dual to one another. This implies, by definition, that the
ambient space and the hypersurface are given in terms of a dual pair of
reflexive polyhedra, with $\D\subset M$ and $\D^*\subset N=\Hom(M,\Z)$ being
exchanged under the mirror involution.
Mirror symmetry thus derives from an elementary combinatorial duality.

Because of the large number of hypersurfaces that exist in these spaces
\cite{c3d,c4d,crp,ams}
only little work was directed towards complete intersections:
A list of transversal configurations for codimension two Calabi--Yau manifolds
in $W\P$ spaces was produced by Klemm \cite{AKcy}. As in the case of
hypersurfaces, there is, however, in general no mirror construction available
in that context \cite{ciwps}. In the toric setup, the mirror construction for
hypersurfaces could be extended to general complete intersections by
Batyrev and Borisov \cite{bb1,bb2}. In addition to a reflexive polyhedron
$\D^*$ that describes the ambient space this involves a decomposition of
$\D$ into a Minkowski sum of polytopes $\D_i$ that are related to the
equations defining the complete intersection. The Calabi--Yau condition
implies that these $\D_i$ are dual to a partition of the vertices of $\D^*$,
which is called nef because the corresponding divisors are numerically
effective. Nef partitions of refelexive polyhedra again feature a beautiful
combinatorial duality that implements the mirror involution, as has been
proven on the level of Hodge data in \cite{strh}.

In the present paper we work out a number of examples of toric
complete intersection Calabi--Yau manifolds and discuss the relation of this
construction to $W\P$ spaces. Identifying CICYs in $W\P$ spaces as a special
case of the toric construction will provide, among other benefits, the mirrors
for these manifolds. In the case of hypersurfaces in $W\IP^4$, the Newton
polytope of a transversal quasihomogeneous polynomial \cite{fl89,cqf}
can be identified with the polyhedron $\D$, whose dual provides the toric
resolution of the ambient space. It is thus clear that, for codimensions
$r>1$, we should look for the identification by trying to relate the
Newton polytopes of the defining polynomial equations of degrees $d_i$
to a nef partition $\D_i$ of some reflexive polyhedron $\D$.

This indeed works for many cases, but the situation is not so straightforward.
Already in the case of hypersurfaces reflexivity of the Newton polytope is
only guaranteed for dimensions $n\le4$ \cite{hs} and indeed
breaks down for the case of Calabi--Yau 4-folds \cite{klry,fft}.
The most relevant situation from the string theory point of view is that
of 3-folds, where already for codimension 2 the Newton polyhedra have
dimension 5. Indeed, already in the second example in the list of ref.
\cite{AKcy}, namely degree (3,4) equations in $\IP_{1,1,1,1,1,2}$, the
Newton polyhedron $\D(7)$ for a degree $7$ equation, is not reflexive. It is,
however, possible to reduce $\D(7)$ to a reflexive polyhedron $\D$ by omitting
5 points, so that its dual provides a toric resolution of singularities of
the weighted projective space. Moreover, in this example, for one of the nef
partitions of $\D$ the Hodge data agree with the $W\P$ result of
\cite{AKcy}.

In general, $\D(d_1+d_2)$ may differ from the Minkowski sum $\D(d_1)+\D(d_2)$
and neither of the two polytopes has to be reflexive.
In many, but not all cases, we can nevertheless find a simple modification of
these polytopes that makes the Hodge data agree, and with some more work one
can check the identification in more detail. In the present note we analyze
a number of examples from the list in \cite{AKcy} and discuss the different
situations that can occur. Apart from our interest in this specific class of
examples, we wrote a program that generates all nef partitions with
codimension two for arbitrary
reflexive polyhedra and that computes the Hodge data for the resulting
Calabi--Yau manifolds. Using the list of 4-fold polyhedra that were obtained
in \cite{fft} we produce a sizeable list of Hodge numbers and compare them
with the complete lists for toric hypersurface. 
Most of the new Hodge numbers lie near the lower ``boundary region'' at
$h^{11}=1$ and appear from a starting polyhedron in the $N-$lattice with
less than 20 points. In particular, we doubled the number of known spectra
with $h^{11}=1$. Altogether we found 87 pairs of new Hodge numbers not
contained in the complete list of toric hypersurfaces \cite{c3d}.
They are listed in table 2 and discussed in section 6.2.

The paper is organized as follows: In section 2 we recall some facts about
toric geometry, mainly to set up our notation. We will use the approach of
the homogeneous coordinate ring, as introduced by Cox \cite{Cox1}. Some basic
facts on the combinatorial data of nef partitions \cite{bor} and a new
criterion for the nef property (proposition \ref{equinef}), which was
used in the computations,
can be found in section 3. In section 4 we recall how to compute a Gorenstein
cone from a nef partition \cite{strh,bb1} and
secton 5 summarizes the polynomials defined in \cite{strh} to calculate the
Hodge data for a Gorenstein cone arising from a nef partition.
The formula used in our program can be found in remark \ref{endformula}.
In section 6 we discuss a number of examples of complete intersections in
$W\IP^5$. We conclude with a discussion of our results, which will be posted
at our web site \cite{KScy} and some of which are listed in the appendix.
A reader who is only interested in new results can take a look at
proposition \ref{equinef} for a new criterion of a nef partition,
section 6.1 for comparing our results with codimension two Calabi--Yau
manifolds in $W\P^5$ spaces \cite{AKcy}, and section 6.2 for new Hodge numbers.

\section{Toric geometry and complete intersections}

Toric geometry is a generalization of projective geometry where the gluing data
of an algebraic variety are encoded in a fan $\S$ of convex rational cones.
Often, the fan is given in terms of (the cones over the faces of) a
polytope $\bar\D$ whose vertices lie on some lattice $N$ \cite{Ful,Oda}.
A very useful way of defining these spaces is to introduce homogeneous
coordinates $z_i$ for all generators $v_i\in N$ $(i=1,\dots,n)$ of the
one-dimensional cones
in $\S$ (e.g. the vertices of $\bar\D$) and to consider the quotient
of $\IC^{n}-Z$ by identifications
$$
        (z_1,\ldots,z_n) \sim (\l^{q^{(I)}_1}z_1,\ldots,\l^{q^{(I)}_n}z_n),
        ~~~~~\sum q^{(I)}_iv_i=0, ~~~I=1,\ldots,n-d,
$$
where the scaling weights $q^{(I)}_i$ describe all linear relations among
the generators $v_i$ and $d=\dim (N)$ is the dimension of the resulting
toric variety $\IP_\S$ \cite{Cox1,cox,ukr}.
In the special case $n=d+1$ of a weighted
projective space the exceptional set $Z$, which is determined in terms of
the fan, only consists of the origin $z_i=0$.

Ample line bundles on $\IP_\S$ correspond to
polytopes $\D$ in the dual space $M=\mao{Hom}(\IZ,N)$ \cite{Ful}.
Toric varieties found their way into string theory when Batyrev \cite{bat}
showed that
the generic section of the line bundle corresponding to $\D$ defines a
Calabi--Yau hypersurface in $\IP_\S$ if $\bar\D$ is equal to the dual
$$
        \D^*=\{x\in N_\IR~|~\langle y,x\rangle\ge-1
                ~\forall y\in\D\subset M_\IR\}
$$
of $\D$, where $N_\IR$ is the real extension of the lattice $N$.
A lattice polytope $\D$ whose dual $\D^*$ is also a lattice polytope is called
reflexive. A necessary condition for this is that the origin is the unique
interior lattice point of $\D$. Moreover, it turned out that
the family of CY hypersurfaces in $\IP_{\S(\D)}$ that is defined by $\D$
provides the mirror family to the family of CY varieties that are based on
$\bar\D=\D^*$  in the sense that the Hodge numbers $h^{p,q}$ and
$h^{d-p,q}$ are exchanged \cite{bat}.
At that time it had just become clear
that hypersurfaces in weighted projective spaces are close to but not exactly
mirror symmetric \cite{nms,KlS}.
This is true even if orbifolds and discrete torsion
        are included \cite{aas,dt}, which do help in the situation where the
        Berglund-Hübsch \cite{be93,mmi} construction applies \cite{ade}.
Beyond the construction of
the missing mirror manifolds, however, Batyrev's results introduced to the
physicist's community beautiful and extremely useful new techniques,
which later turned out also to apply to the analysis of fibration
structures that are important in string dualities \cite{AL,Fv,Arev,cafo}:
In toric geometry CY fibrations manifest themselves as reflexive sections or
projections of the polytopes $\D^*$ and $\D$, respectively \cite{k3,fft}.

\section{Nef partitions}

In the case of a hypersurface, the supporting polyhedron of the generic section
of an ample line bundle on $\P_\Sigma$ must be reflexive in order to get a
Calabi--Yau hypersurface in $\P_\Sigma$. To generalize this condition to the
case of codimension $r \geq 1 $, i.e. to ensure that the intersection of
$r$ hypersurfaces is a Calabi--Yau manifold, the reflexive polytope $\Delta
\subset M_\R$ must fulfil the so called nef condition \cite{bor}. In this
section, we will shortly discuss the combinatorial properties of nef partitions
and give a new criterion for a reflexive polytope to decompose into a
nef partition, proposition \ref{equinef}, which can be used to calculate
these partitions in a simple way, as described in remark \ref{calcnef}.

Let $\Delta \subset M_\R$ be a reflexive polytope and $\Delta^\ast \subset
N_\R$ it's dual. From now on we denote by $\Delta^v$ the {\em set of vertices}
of a polytope $\Delta$. Let $E := {\Delta^\ast}^v$ be the set of vertices of
$\Delta^\ast$. We define the $d-$dimensional complete {\em fan}
$\Sigma[\Delta^\ast]$ as the union of the zero-dimensional cone $\{0\}$
together with the set of all cones
$$C[F]=\{0\}\cup\{z\in N_\R:\lambda z\in F\ \mathrm{for}\ \mathrm{some}\
\lambda\in\R_>\}$$
that support faces $F$ of $\Delta^\ast$. Assume that
there exists a representation of
$E=E_1\cup\dots\cup E_r$ as the union of disjoint subsets $E_1,\dots,E_r$ and
integral convex $\Sigma[\Delta^\ast]-$ piecewise linear support functions
$\varphi_i:N_\R\rightarrow \R\ (i=1,\dots,r)$  such that
$$
\varphi_i(e)=
\b{cases}
1& \text{if $e\in E_i$,}\\
0& \text{otherwise.}
\e{cases}
$$
Each $\varphi_i$ corresponds to a line bundle $L_i$
that defines a supporting polyhedron $\Delta_i$ for the global sections:
$$
\Delta_i = \{\d z\in M_\R : \langle \d z, z\rangle \ge -\varphi_i(z)\ \forall \
z\in N_\R\}.
$$
        $\varphi_i$ defines $L_i$ in the following way:
For each cone of
        maximal dimension $C$ there is a $m_C\in M$ such that
        $\varphi_i\mid_C = m_C\mid_C$, where the $m_C$ have to coincide on the
        intersection of two cones. Since the fan is complete, $\varphi_i$ is
        determined uniquely by the set $\{m_C\}$. The
        line bundle $L_i$ is then given by the data $(U_C,\chi(m_C))$, where
        $\{U_C\}$ is an open covering of the toric variety with open sets
        corresponding to the cones of maximal dimension and $\chi(m_C)$ can
        be regarded as a monomial $x^{m_C}$. The important point is that the
        transition functions $\chi(m_{\tilde C}-m_C)$ arising from this
        construction are regular on the intersection of the corresponding
        open sets (for details see \cite{Oda} or \cite{Ful}).

Conversely, each function $\varphi_i$ is uniquely defined by the polyhedron
$\Delta_i$.
A {\em Calabi--Yau complete intersection} is then determined  by the
intersection of the closure of $r$ hypersurfaces, each corresponding to a
global section of a line bundle $L_i$ \cite{strh,bb1}.

\b{Def} \label{nef}
If there exists a reflexive polytope $\Delta$ and $r$ functions
$\varphi_1,\dots,\varphi_r$ as defined above, we call the data
$$
\Pi(\Delta)=\{\Delta_1,\dots,\Delta_r\}
$$
a {\em nef partition}.
\e{Def}
Equivalent to $\Pi(\Delta):=\{\Delta_1,\dots,\Delta_r\}$ being a nef partition
is that any two $\Delta_i$ only have $\{0\}$ as a common point and that
$\Delta$ can
be written as the Minkowski sum $\Delta_1+\dots+\Delta_r = \Delta$, as is
shown by the following proposition:
\b{Pro}\label{equinef}
$\Pi(\Delta)=\{\Delta_1,\dots,\Delta_r\}$ is a nef partition if and only if
$\Delta$ is the Minkowski sum of $r$ rational polyhedra
$\Delta=\Delta_1+\dots+\Delta_r$ and $\Delta_i\cap\Delta_j=\{0\}\ \forall\
i\neq j$
\e{Pro}

\noindent \pro\\
$\Rightarrow$: Assume that $\Delta$ can be written as the Minkowski sum of $r$
rational polyhedra $\Delta=\Delta_1+\dots+\Delta_r$ with
$\Delta_i\cap\Delta_j=\{0\}\ \forall\ i\neq j$. Define $r$ functions
$\varphi_i:N_\R\rightarrow \R$ as
$$\varphi_i(z)=-\min\limits_{\d z\in \Delta_i}\langle\d z, z\rangle\ \forall\
z\in N_\R.$$
\b{itemize}
   \item
The $\varphi_i$ are linear on cones of $\Sigma[\Delta^\ast]$: It is sufficient
to consider restrictions of the $\varphi_i$ to cones of maximal dimension
$C[F]$, where
$$F=\Delta^\ast\cap\{z\in N_\R\ :\langle \d e,z\rangle=-1\}$$
is a facet of $\Delta^\ast$ corresponding to a vertex $\d e\in\Delta^v$. Now
let $\d e=\d e_1+\dots+\d e_i+\dots+\d e_r$, where $\d e_i\in \Delta_i^v $
denotes a  vertex of $\Delta_i\ (i=1,\dots,r)$. If we take another vertex ${\d
e_i}^\prime \neq \d e_i\in \Delta_i^v$, then the sum ${\d e}^\prime=\d
e_1+\dots+{\d e_i}^\prime +\dots+\d e_r$ denotes another vertex of $\Delta$.
Clearly, $\langle\d e, z\rangle\le \langle{\d e}^\prime, z\rangle\ \forall z\in
C[F]$, i.e. $\ \langle\d e_i, z\rangle\le \langle{\d e_i}^\prime, z\rangle\
\forall z\in C[F]$. Hence $\varphi_i(z)=-\langle\d e_i,z\rangle\ \forall\ z\in
C[F]$.
    \item
Convexity of all $\varphi_i$ follows immediately from their definition.
   \item
$\varphi_i(e) \in \{0,1\}\ \forall\  e\in E,\ i=1,\dots, r$: For every
function
$\varphi_i$ we observe that $0\in \Delta_i$ implies $\varphi_i \ge 0$ and
$\Delta_i\subseteq\Delta$ implies $\varphi_i(e) \le 1\ \forall\  e\in E$.
   \item
$\varphi_i(e)=1\Rightarrow \varphi_j(e)=0\ \forall\ j\neq i$: Assume
$\varphi_i(e)=\varphi_j(e)=1$ for $i\neq j$ $\Rightarrow\exists \d z_i\in
\Delta_i,\ \d z_j\in \Delta_j: \langle\d z_i, e\rangle=\langle\d z_j,
e\rangle=-1\Rightarrow\exists \d z=\d z_i+\d z_j\in \Delta$ with $\langle\d
z,e\rangle=-2$. This contradicts $\Delta^\ast$ being dual to $\Delta$.
   \item
$\forall e\in E\ \exists\ i\in\{1,\dots,r\}$ with $\varphi_i(e)=1:$ Assume
$\exists e\in E:\varphi_i(e)=0\ \forall\ i=1,\dots,r.$ By duality of $\Delta$
and $\Delta^\ast\ \exists\d z\in \Delta:\langle\d z,e\rangle=-1$, where $\d z$
is contained in the facet dual to $e$. Now $\d z=\d z_1+\dots+\d z_r$ with $\d
z_i\in \Delta_i\ \forall\ i=1,\dots,r$. $\Rightarrow\ \exists \d z_k\in
\Delta_k$ with $\langle\d z_k,e\rangle<0.$ This  contradicts $\varphi_i(e)=0\
\forall\ i=1,\dots,r.$
\e{itemize}
$\Leftarrow$: Follows from $$\Delta=\{\d z\in M_\R:\langle\d
z,z\rangle\geq-\varphi(z)\ \forall\ z\in N_\R\},$$ where
$\varphi=\sum_i\varphi_i\ (i=1,\dots, r)$.\hfill $\Box$\\

It can be shown that every nef partiton of a reflexive polytope $\Delta$ gives
a dual nef partition of a reflexive polytope $\nabla$, which turns out to be an
involution on the set of  reflexive polytopes with nef partitions:

\b{Rem}\cite{bor}\label{dnef}
Let $\ \Pi(\Delta)=\{\Delta_1,\dots,\Delta_r\}$ be a nef partition and denote
by $E=E_1\cup\dots\cup E_r$ the set of vertices ${\Delta^\ast}^v$.
Define $r$ rational polyhedra $\nabla_i\subset N_\R\ (i=1,\dots,r)$ as
$$\nabla_i=\mathrm{Conv}(E_i\cup\{0\}).$$
Then there is the following relation between $\Delta_i$ and $\nabla_j\ (i,j =
1,\dots,r)$:
$$
\langle \Delta_i,\nabla_j\rangle=
\b{cases}
\ge -1& \text{if}\ \ i=j\\
\ge 0& \text{otherwise,}
\e{cases}
$$
and the $\nabla_i$ are maximal with that property. In particular
$\nabla=\nabla_1+\dots+\nabla_r$ is a reflexive polyhedron with a nef partition
$\ \Pi(\nabla)=\{\nabla_1,\dots,\nabla_r\}$, and there is a natural involution
on the set of reflexive polyhedra with nef partitions:
$$
\iota\ :\  \Pi(\Delta)=\{\Delta_1,\dots,\Delta_r\}\mapsto
\Pi(\nabla)=\{\nabla_1,\dots,\nabla_r\}.
$$
\e{Rem}

\b{Rem}\label{calcnef}
The following procedure can be used to find all nef partions of a reflexive
polyhedron $\Delta\subset M_\R$:
\begin{itemize}
    \item
First calculate $\Delta^\ast\subset N_\R$.
    \item
Take disjoint unions $E=E_1\cup\dots\cup E_r$ of vertices of $\Delta^\ast$.
   \item
Check if $\nabla=\nabla_1+\dots+\nabla_r$ with $\nabla_i=\Co(E_i\cup \{0\})$ 
is reflexive and $\nabla_i\cap \nabla_j = \{0\}\ \forall\ i\neq j$.
\end{itemize}
\e{Rem}

\section{Gorenstein cones}

The (string theoretic) Hodge numbers of a 
Calabi--Yau manifold corresponding to a nef partition
are the coefficients of the  $E-$polynomial
\[      E_{st}(V;u,v)=\sum (-1)^{p+q} h_{st}^{p,q}u^pv^q,
\]
which can be computed from a higher-dimensional Gorenstein cone \cite{strh} 
that is constructed using the data of a nef partition \cite{bb1,strh}.
In this section we will give the definition of a Gorenstein cone 
and recall its construction starting with a nef partition.

A rational cone $C\subset M_\R$ is called {\em Gorenstein} if there exists a
point $n\in N$ in the dual lattice such that $\langle v, n\rangle=1$ for all
generators of the semigroup $C\cap M$. Given a nef partition $\
\Pi(\Delta)=\{\Delta_1,\dots,\Delta_r\}$, we can construct such a cone. First
we go to a larger space and extend the canonical pairing: Let $\Z^r$ be the
standard $r-$dimensional lattice and $\R^r$ its real scalar extension. We put
$\bar N=\Z^r\oplus N,\ \bar d = d+r$ and $\bar M = \mathrm{Hom}(\bar N,\Z)$. We
extend the canonical $\Z-$bilinear pairing $\langle \ast ,\ast\rangle :M \times
N \to\Z$ to a pairing between $\bar M$ and $\bar N=\Z^r\oplus N$ by the formula
$$
\langle(a_1,\dots,a_r,m),(b_1,\dots,b_r,n)\rangle=\sum\limits_{i=1}^{r}a_i
b_i+\langle m,n\rangle.\label{ext}
$$
The real scalar extensions of $\bar N$ and $\bar M$ are denoted by $\bar
N_\R$ and $\bar M_\R$, respectively, 
with the corresponding $\R-$bilinear pairing
$\langle\ast,\ast\rangle: \bar M_\R\times \bar N_\R\to \R$.

\b{Def}\label{gorenef}
{}For a nef partiton $\ \Pi(\Delta)=\{\Delta_1,\dots,\Delta_r\}$ we construct
a $\bar d-$dimensional Gorenstein cone $C_\Delta \subset \bar M_\R$
$$
C_\Delta =\{(\lambda_1,\dots ,\lambda_r,\lambda_1 \d z_1+\dots +\lambda_r \d
z_r)\in\bar M_\R:\lambda_i\in \R_\ge, \d z_i\in \Delta_i,\ i=1,\dots r\},
$$
with $n_\Delta\in \bar N$ uniquely defined by the conditions
\bea
\langle \d z,n_\Delta\rangle &=& 0\ \forall\ \d z\in M_\R\subset\bar
M_\R\nonumber\\
\langle \d e_i,n_\Delta\rangle &=& 1\ \mathrm{for}\ i=1,\dots,r,\nonumber
\eea
where $\{\d e_1,\dots,\d e_r\}$ is the standard basis of $\Z^r\subset \bar M$.
\e{Def}

Note that all generators of $C_\Delta\cap \bar M$ lie on the hyperplane 
$\langle \d z,n_\Delta \rangle = 1$. They span the $\bar d -1-$dimensional
supporting polyhedron
 $$ K_\Delta=\{\d z\in C_\Delta:\langle\d z,n_\Delta\rangle=1\}$$
of $C$.
Since $K_\Delta\cap\bar M$ has no interior point, we get
$$K_\Delta\cap\bar M=\bigcup\limits_{i=1,\dots,r}(\d e_i\times
\Delta_i)\cap\bar M.$$

\b{Rem}\cite{bb1}\label{dgorenef}
Let $\ \Pi(\nabla)=\{\nabla_1,\dots,\nabla_r\}$ be the dual nef partition. Then
the Gorenstein cone
$$
C_\nabla=\{(\mu_1,\dots ,\mu_r,\mu_1 z_1+\dots +\mu_r z_r)\in\bar N_\R:\mu_i\in
\R_\ge, z_i\in \nabla_i,\ i=1,\dots r\}
$$
is dual to $C_\Delta$ defined in \ref{gorenef}. Note, however, that 
$K_\Delta$ is {\em not} dual to $K_\nabla$!
\e{Rem}

\section{Combinatorical polynomials of Eulerian posets}

Batyrev and Borisov gave an explicit formula for the string-theoretic
$E-$polynomial for a Calabi--Yau complete intersection $V$ in a Gorenstein
toric Fano variety \cite{strh}. This polynomial depends only on the
combinatorial data of the corresponding Gorenstein cone. We will
give some basic definitions of combinatorial polynomials on Eulerian Posets,
which are used to compute the $E-$polynomial, and formulate it in a way which 
can be used for the explicit calculation of the Hodge numbers.

Let $P$ be an {\em Eulerian Poset}, i.e. a finite partially ordered set with
unique minimal element $\hat 0$, maximal element $\hat 1$ and the
same length $d$
of every maximal chain of $P$. For any $x\leq y\in P$, define the {\em
interval} $I=[x,y]$ as
$$[x,y]=\{z\in P:x\leq z\leq y\}.$$
In particular, we have $P=[\hat 0,\hat 1]$. Define the {\em rank function}
$\rho:P\rightarrow\{0,\dots,d\}$ on $P$ by setting $\rho(x)$ equal to the
length of the interval $[\hat 0,x]$. Note that for any Eulerian Poset $P$,
every interval $I=[x,y]$ is again an Eulerian Poset with rank function
$\rho(z)-\rho(x)\ \forall z\in I$. If an Eulerian Poset has rank $d$, then the
{\em dual Poset} $P^\ast$ is also an Eulerian Poset with rank function
$\rho^\ast=d-\rho$.

\b{Exa}\label{Cpos}
Let $C\in N_\R$ be a $d-$dimensional cone with its dual $C^\ast\in M_\R$.
There is a canonical bijective correspondence $F\leftrightarrow F^\ast$ between
faces $F\subseteq C$ and $F^\ast\subseteq C^\ast\ $ with
$\mathrm{dim}F+\mathrm{dim}F^\ast=d$ \cite{Ful},
$$F\mapsto F^\ast=\{z\in C^\ast:\langle\d z,z\rangle=0\ \forall z\in F\ \},$$
which reverses the inclusion relation between faces. We denote the faces of $C$
by indices $x$ and define the poset $P=[\hat 0,\hat 1]$ as the poset of all
faces $C_x\subseteq C$ with maximal element $C$ and minimal element $\{0\}$ and
rank function $\rho(x)=\mathrm{dim}(C_x)\ \forall x\in P$. The dual poset
$P^\ast$ can be identified with the poset of faces $C^\ast_x\subseteq C^\ast$
of the dual cone $C^\ast$ with rank function
$\rho^\ast(x^\ast)=\mathrm{dim}(C^\ast_x)\ \forall\ x^\ast\in P^\ast$.
\e{Exa}

\b{Def}\label{B}
Let $P$ be an Eulerian Poset of rank $d$ as above. Define the {\em polynomial}
$B(P;u,v)\in \Z[u,v]$ by the following rules \cite{stan,strh}:
\b{itemize}
   \item
   $B(P;u,v) = 1$ if $d=0$;
   \item
   the degree of $B(P;u,v)$ with respect to $v$ is less than $d/2$;
   \item
   $\sum\limits_{\hat 0\leq x\leq\hat 1}\!\!\! B([\hat
0,x];u^{-1},v^{-1})(uv)^{\rho (x)}(v-u)^{d-\rho (x)}=\!\sum\limits_{\hat 0\leq
x\leq\hat 1}\!\!\! B([x,\hat 1];u,v)(uv-1)^{\rho (x)}.$
\e{itemize}
\e{Def}

Let us consider how we can construct the $B-$polynomial for an interval
$I=[x,y]\subseteq P$ with $d=\rho(y)-\rho(x)$:\
Suppose we know the $B-$polynomials $B(\tilde I;u,v)$ for all sub-intervals
$\tilde I=[\tilde x,\tilde y]\subset I$.
Then we know all terms of the relation formula for the $B-$polynomials in
\ref{B} except for 
$B(I;u,v)$ on the right hand side and $B(I;u^{-1},v^{-1})(uv)^d$
on the left hand side. Because the $v-$degree of $B(I;u,v)$ is less than $d/2$,
the possible degrees of monomials with respect to $v$ in $B(I;u,v)$ and
$B(I;u^{-1},v^{-1})(uv)^d$ do not coincide and we can calculate $B(I;u,v)$. So
if we have to compute $B(P;u,v)$, we first have to calculate the
$B-$polynomials for all intervals with rank $0$ (which are per definition 1),
then those intervals with rank $1$, and so on.

\b{Rem}\label{Bd}
Let $P$ be an Eulerian Poset of rank $d$, $P^\ast$ be the dual. Then the
polynomial defined in (\ref{B}) satisfies
$$B(P;u,v)=(-u)^d B(P^\ast;u^{-1},v).$$
\e{Rem}

\b{Def}\label{S}
Let $P$ be the Eulerian Poset corresponding to the Gorenstein cone $C=C_\Delta
\subset \bar M_\R$ from definition \ref{gorenef}. Define two functions on the
set of faces of $C$ by
\bea
S({C}_x,t) &=& (1-t)^{\rho(x)}\sum\limits_{m\in C_x\cap \bar
M}t^{\mathrm{deg}(m)}\nonumber\\
T(C_x,t) &=& (1-t)^{\rho(x)}\sum\limits_{m\in \mathrm{Int}(C_x)\cap \bar
M}t^{\mathrm{deg}(m)},\nonumber
\eea
where $\mathrm{Int}(C_x)$ denotes the relative interior of $C_x\subseteq C$ and
$\mathrm{deg}(m)=\langle m,n_\Delta \rangle$.
\e{Def}

The following statement is a consequence of the Serre duality \cite{dako}:

\b{Pro}\label{Serre}
For the Gorenstein cone $C=C_\Delta \subset \bar M_\R$ the functons $S$ and $T$
are polynomials: $S(C_x,t),\ T(C_x,t)\in \Z[t]$, and they satisfy the relation
$$S(C_x,t)=t^{\rho(x)}T(C_x,t^{-1}).$$
\e{Pro}

\b{Rem}\label{rel}
For $S=\sum_i a_i t^i$ and $T=\sum_i b_i t^i$ as defined above \ref{Serre}
implies that
$$a_0+a_1t+\dots+a_nt^n=b_0t^n+b_1t^{n-1}+\dots+b_{n-1}t+b_n,$$
where $n=\di C_x$ and we get the relations
$$a_i=b_{n-i}\ (i=1,\dots,n)$$
for the coefficients of $S$ and $T$. Since $a_0=1$ and $b_0=0$, the leading
coefficients are determined to be $a_n=0$ and $b_n=1$. So it is sufficient to
calculate $\mid C_x\cap m\cdot K_\Delta\mid$ and $\mid \Int(C_x\cap m\cdot
K_\Delta)\mid$ for $m=0,\dots ,[\dim(C_x)/2]$ and to use the fact that
$a_i=b_{n-i}$ for $i > \dim(C_x)/2$.
\e{Rem}


Batyrev and Borisov showed in their paper \cite{strh} that the {\em
string-theoretic E-polynomial} of a nef partition can be calculated from the
data of the corresponding Gorenstein cone:

\b{Pro}\label{E}
\cite{strh} Let $\ \Pi(\Delta)=\{\Delta_1,\dots,\Delta_r\}$ be a nef partition
and $C=C_\Delta\subset \bar M_\R$ be
the $\bar d-$dimensional reflexive Gorenstein
cone defined in \ref{gorenef} (with dual cone $C^\ast = C_\nabla\subset \bar
N_\R$). Denote by $P$ the poset of faces $C_x\subseteq C$ (see example
\ref{Cpos}). Then the string-theoretic $E-$polynomial is given by

$$
E_{st}(V;u,v) = \sum\limits_{I=[x,y]\subseteq P}\frac{(-1)^{\rho
(y)}}{(uv)^r}(v-u)^{\rho (x)} B(I^\ast;u,v)(uv-1)^{\bar d-\rho
(y)}A_{(x,y)}(u,v),
$$
with
$$
A_{(x,y)}(u,v) = \sum\limits_{(m,n)\in \mathrm{Int}(C_x)\cap \bar
M\times\mathrm{Int}(C_y^\ast)\cap \bar N}\left(\frac{u}{v}\right)^{\mathrm{deg}
(m)}\left(\frac{1}{uv}\right)^{\mathrm{deg} (n)}.
$$
\e{Pro}

The dual partition $\ \Pi(\nabla)=\{\nabla_1,\dots,\nabla_r\}$ corresponds to
the Calabi--Yau complete intersection $W$ and $(V,W)$ is a {\em mirror
pair} of (singular) Calabi--Yau varieties, at least in the sense that
$$
 E_{st}(V;u,v) = (-u)^{d-r}E_{st}(W;\frac{1}{u},v),
$$
or equivalently, $h_{st}^{p,q}(V)=h_{st}^{n-p,q}(W)$ for 
$0\leq p,q\leq n=\dim(V)=\dim(W)$.

\b{Rem}\label{endformula}
Using the duality \ref{Bd} for the $B-$polynomials and definition \ref{S} with
relation \ref{Serre} between the $S-$ and $T-$polynomials, we can write the
$E-$polynomial as
$$
E_{st}(V;u,v) = \sum\limits_{I=[x,y]\subseteq P}\frac{(-1)^{\rho (x)}u^{\rho
(y)}}{(uv)^r}S\left( C_x,\frac{v}{u}\right) S\left( C_y^\ast,uv\right)
B(I;u^{-1},v).
$$
This equation can be used for explicit calculations.
\e{Rem}

 \newpage

\section{Results}

Using the formula for the $E$-polynomial \ref{endformula} we are now able to
construct Calabi--Yau complete intersections starting with a reflexive polytope
$\Delta\subset M_\R$ (or $\Delta^\ast\subset N_\R$). Our first task is to
compare the toric construction to a list of complete intersections in 
weighted projective spaces, which was produced 
by Klemm \cite{AKcy}. Then we construct a large number of
nef partitions for different classes of five-dimensional reflexive polytopes 
and compare the Hodge data with the complete restults that are available for
toric hypersurfaces \cite{c4d,KScy}.

\subsection{Comparison with Weighted Projective Space}

In order to identify complete
intersections in $W\P^5$ as special cases of the toric construction it is
natural to start with the Newton polyhedron and
to compare the Hodge data for various nef partitions.
In what follows we will analyze some
examples from Klemm's list \cite{AKcy} and discuss the different situations 
that can occur.

In the simplest case the Newton polyhedron $\D(d)$ 
corresponding to degree
$(d_1,d_2)$ equations with $d=d_1+d_2$ 
is reflexive and the Hodge numbers of a nef partition
$\Pi(\Delta(d))=\{\tilde \Delta_1,\tilde\Delta_2\}$ agree with those given in
\cite{AKcy}. This works already for the first example of degree $(4,2)$
equations in $\P_{1,1,1,1,1,1}$.
In general, $\D(d_1+d_2)$ may differ from the Minkowski sum $\D(d_1)+\D(d_2)$
and none of the two polytopes has to be reflexive. In many cases we find a
simple modification of these polytopes that makes the Hodge data agree:

\begin{itemize}
\item
Already in the second example of this list the Newton polyhedron $\D(7)$ for
the weight system of $W\P_{1,1,1,1,1,2}$ 
is not reflexive. It is, however,
possible to reduce $\D(7)$ to a reflexive polyhedron $\D$ by omitting 5
points, so that its dual provides a toric resolution of singularities of the
weighted projective space. Indeed, the Hodge data for $(d_1,d_2)=(3,4)$ 
matches for one nef partition of the resulting polytope.
\item
Another possibility is that the Newton polyhedron is reflexive, but the Hodge
numbers do not agree. In such a case we can compute the Minkowski sum
$\tilde \D= \D(d_1)+\D(d_2)$ and check if it is reflexive and gives the right
Hodge numbers. This works, for example, for degrees $(d_1,d_2)=(5,3)$ in case
of the weight system for $W\P_{1,1,1,1,2,2}$.
\end{itemize}

There are still some examples where we were not able to reproduce the Hodge
data. For example, in case of the weight system for $W\P_{1,1,1,1,2,3}$
and degrees $(d_1,d_2)=(5,4)$, neither
$\Delta(9)$ (which has 575 points) nor the Minkowski sum $\D(5)+\D(4)$ 
(with only 211 points) is reflexive. The largest reflexive subpolytope of 
$\Delta(9)$ has 570 points, but it's nef partitions do not yield the right
Hodge numbers. Omitting up to 30 points we find another 21 reflexive polyhedra,
but non of their nef partitions yields $h^{11}=2$ and $h^{12}=84$. We thus
found no candidate for a toric description and a more detailed analysis of
the geometry would be required to check if a toric description exists.

\subsection{New Hodge numbers}

Of course, one of our main interests is to find new Hodge data. 
In \cite{c4d,KScy} the complete set of 30108 pairs of Hodge numbers 
corresponding to hypersurfaces in toric varieties has been constructed
(the number is 15122 if we count those with $h^{11}\leq h^{12}$ because 
136 are selfdual). 
Picking out only the new Hodge data from the list of 2387
pairs arising from weighted $\P^5$ \cite{AKcy}, i.e. those not contained in 
\cite{KScy}, there remain only 15 new data, which we list in table 1. 
(Note that this class is not mirror symmetric.)

\begin{figure}\begin{center}
\b{tabular}{|r|r|r||r|r|r||r|r|r|}
\hline
$h^{11}$ & $h^{21}$ &R& $h^{11}$ & $h^{21}$ &R& $h^{11}$ & $h^{21}$&R
\\\hline\hline
1 & 61 & x & 2 & 62 & x & 7 &  26 &  \\\hline
1 & 73 & x & 2 & 68 & x & 8 &  20 &  \\\hline
1 & 79 & x & 3 & 47 & x &12 &  12 &  \\\hline
1 & 89 & x & 3 & 55 & x &13 &  13 & x\\\hline
1 &129 & x & 3 & 61 & x &17 &  11 &  \\\hline
\e{tabular}
\\[5mm]
Table 1: New Hodge numbers in \cite{AKcy}, as compared to toric hypersurfaces.
\\\hspace*{19truemm} R=x means that we found the same Hodge data
for nef partitions.
\end{center}\end{figure}

Using the toric construction, we started with reflexive polyhedra that
are described by single or combined weight systems, as they were
constructed systematically for Calabi--Yau fourfolds \cite{fft,KScy,ams}, and
from Minkowski sums of Newton polytopes that arise in the context of
weighted projective spaces. In this way we found 16 (with mirror duality 32) 
pairs of new Hodge numbers. They are listed in appendix A, tables 3-5, 
together 
with a detailed information about the starting polyhedron. Most of them lie 
in the lower boundary region $h^{11} \leq 6$, which is less covered by the
``background'' of toric hypersufaces. It is remarkable that almost every pair
of new Hodge numbers corresponds to a starting polyhedron $\Delta^\ast\cap N$
in the $N$-lattice with less than $20$ points. Thus, to get a more complete
result for new spectra in that range, we used the program package that was
written for the classification of reflexive polyhedra \cite{KScy,palp}
to construct a fairly complete set of reflexive polyhedra with up to 10 
points (they were all found as
subpolytopes of some 10000 polyhedra with up to 40 points originating from 
transveral weight systems \cite{fft}).
Indeed, using these polytopes for $\D^*\subset N$, we found $87$ pairs of 
Hodge numbers not contained in \cite{KScy}. They are listed in table $2$ and 
are shown in figure 1
in the background of toric hypersurfaces and CICYs in $W\P^5$ in the 
range of $1\leq h^{11}\leq 10$ and $1\leq h^{21}\leq 170$.

\begin{center}
\begin{tabular}{|r|r||r|r||r|r||r|r||r|r||r|r||r|r|}
\hline
$h^{11}$ & $h^{21}$ & $h^{11}$ & $h^{21}$ & $h^{11}$ & $h^{21}$ & $h^{11}$ &
$h^{21}$ & $h^{11}$ & $h^{21}$ & $h^{11}$ & $h^{21}$\\\hline\hline
 1 &  25 &  2 &  60 &  3 &  24 &  3 &  52 &  4 &  16 &  4 &  53\\\hline 
 1 &  37 &  2 &  62 &  3 &  27 &  3 &  53 &  4 &  22 &  4 &  75\\\hline 
 1 &  61 &  2 &  64 &  3 &  29 &  3 &  54 &  4 &  24 &  4 &  83\\\hline
 1 &  73 &  2 &  66 &  3 &  31 &  3 &  55 &  4 &  26 &  5 &  25\\\hline 
 1 &  79 &  2 &  68 &  3 &  33 &  3 &  56 &  4 &  30 &  5 &  27\\\hline 
 1 &  89 &  2 &  70 &  3 &  35 &  3 &  58 &  4 &  31 &  5 & 102\\\hline 
 1 & 129 &  2 &  72 &  3 &  37 &  3 &  60 &  4 &  32 &  6 &  20\\\hline 
 2 &  30 &  2 &  76 &  3 &  39 &  3 &  61 &  4 &  33 &  6 &  24\\\hline  
 2 &  36 &  2 &  77 &  3 &  41 &  3 &  62 &  4 &  38 &  7 &  22\\\hline 
 2 &  44 &  2 &  78 &  3 &  42 &  3 &  64 &  4 &  39 &  8 &  14\\\hline 
 2 &  50 &  2 &  80 &  3 &  44 &  3 &  68 &  4 &  41 & 13 &  13\\\hline 
 2 &  54 &  2 &  82 &  3 &  47 &  3 &  70 &  4 &  43 & 13 &  15\\\hline      
 2 &  56 &  2 & 100 &  3 &  48 &  3 &  80 &  4 &  45 &    &\\\hline      
 2 &  58 &  2 & 112 &  3 &  49 &  3 & 101 &  4 &  47 &    &\\\hline      
 2 &  59 &  3 &  23 &  3 &  50 &  3 & 113 &  4 &  51 &    &\\\hline   
\end{tabular}
\\[5mm]
Table 2. New Hodge numbers with toric CICYs.
\end{center}

\centerline{\unitlength=2.7pt\begin{picture}(150,50)(15,-25)
        \newcommand{\DI}{\makebox(-5,0){$\diamond$}}
        \def\i#1.#2.{\put(#2,#1){\circle*{.8}}}   \put(0,0){\vector(1,0){175}}
        \def\h#1.#2.{\put(#2,#1){\circle{.75}}}   \put(0,0){\vector(0,1){10}}
        \def\w#1#2{\put(#2,#1){\makebox(0,0){\put(2.5,-0.2){\bf{\DI}}}}}
        \def\hlab#1{\put(#1,-1){\line(0,1){2}}} \hlab{10}\hlab{20}\hlab{30}
        \hlab{40}\hlab{50}\hlab{60}\hlab{70}\hlab{80}\hlab{90}\hlab{100}
        \hlab{110}\hlab{120}\hlab{130}\hlab{140}\hlab{150}\hlab{160}\hlab{170}
        \put(48.5,-4){\small\bf50}\put(97.5,-4){\small\bf100}
        \put(147.5,-4){\small\bf150}\put(175,-4){$\mathbf{h^{21}}$}
        \put(2,8){$\mathbf{h^{11}}$}
%
%
    \put(26,-14){\begin{picture}(0,0)\put(-18,2.2){Figure 1.}\end{picture}
        \parbox{125mm}{%
        87 new toric CICYs (\makebox(.9,.2)[lb]{\put(.5,.95){\circle*{.8}}})
        and CICYs in $W\P^5$ (\makebox(2.17,0)[lb]{\put(3.6,.75){\DI}})
        in the background of toric hypersurfaces 
        (\makebox(.9,.2)[lb]{\put(.4,.95){\circle{.75}}})
        with $h^{11}\leq 10$ and $h^{21}\leq 170$.}}
\w{1}{61}\w{1}{73}\w{1}{79}\w{1}{89}\w{1}{129}\w{2}{62}\w{2}{68}
\w{3}{47}\w{3}{55}\w{3}{61}\w{7}{26}\w{8}{20}
\i1.25.\i1.37.\i1.61.\i1.73.\i1.79.\i1.89.\i1.129.\i2.30.\i2.36.\i2.44.
\i2.50.\i2.54.\i2.56.\i2.58.\i2.59.\i2.60.\i2.62.\i2.64.\i2.66.\i2.68.
\i2.70.\i2.72.\i2.76.\i2.77.\i2.78.\i2.80.\i2.82.\i2.100.\i2.112.\i3.23.
\i3.24.\i3.27.\i3.29.\i3.31.\i3.33.\i3.35.\i3.37.\i3.39.\i3.41.\i3.42.
\i3.44.\i3.47.\i3.48.\i3.49.\i3.50.\i3.52.\i3.53.\i3.54.\i3.55.\i3.56.
\i3.58.\i3.60.\i3.61.\i3.62.\i3.64.\i3.68.\i3.70.\i3.80.\i3.101.\i3.113.
\i4.16.\i4.22.\i4.24.\i4.26.\i4.30.\i4.31.\i4.32.\i4.33.\i4.38.\i4.39.
\i4.41.\i4.43.\i4.45.\i4.47.\i4.51.\i4.53.\i4.75.\i4.83.\i5.25.\i5.27.
\i5.102.\i6.20.\i6.24.\i7.22.\i8.14.
\h1.21.\h1.101.\h1.103.\h1.145.\h1.149.\h2.29.\h2.38.\h2.74.\h2.83.\h2.84.
\h2.86.\h2.90.\h2.92.\h2.95.\h2.102.\h2.106.\h2.116.\h2.120.\h2.122.\h2.128.
\h2.132.\h2.144.\h2.272.\h3.43.\h3.45.\h3.51.\h3.57.\h3.59.\h3.63.\h3.65.
\h3.66.\h3.67.\h3.69.\h3.71.\h3.72.\h3.73.\h3.75.\h3.76.\h3.77.\h3.78.\h3.79.
\h3.81.\h3.83.\h3.84.\h3.85.\h3.87.\h3.89.\h3.91.\h3.93.\h3.95.\h3.99.\h3.103.
\h3.105.\h3.107.\h3.111.\h3.115.\h3.119.\h3.123.\h3.127.\h3.131.\h3.141.
\h3.165.\h4.28.\h4.34.\h4.36.\h4.37.\h4.40.\h4.42.\h4.44.\h4.46.\h4.48.\h4.49.
\h4.50.\h4.52.\h4.54.\h4.55.\h4.56.\h4.57.\h4.58.\h4.59.\h4.60.\h4.61.\h4.62.
\h4.63.\h4.64.\h4.65.\h4.66.\h4.67.\h4.68.\h4.69.\h4.70.\h4.71.\h4.72.\h4.73.
\h4.74.\h4.76.\h4.77.\h4.78.\h4.79.\h4.80.\h4.81.\h4.82.\h4.84.\h4.85.\h4.86.
\h4.88.\h4.89.\h4.90.\h4.91.\h4.92.\h4.93.\h4.94.\h4.96.\h4.97.\h4.98.\h4.100.
\h4.101.\h4.102.\h4.104.\h4.106.\h4.108.\h4.109.\h4.110.\h4.112.\h4.114.
\h4.116.\h4.118.\h4.120.\h4.121.\h4.122.\h4.124.\h4.126.\h4.128.\h4.130.
\h4.136.\h4.142.\h4.144.\h4.148.\h4.154.\h4.162.\h4.166.\h5.20.\h5.29.\h5.31.
\h5.33.\h5.35.\h5.37.\h5.38.\h5.39.\h5.40.\h5.41.\h5.42.\h5.43.\h5.44.\h5.45.
\h5.46.\h5.47.\h5.48.\h5.49.\h5.50.\h5.51.\h5.52.\h5.53.\h5.54.\h5.55.\h5.56.
\h5.57.\h5.58.\h5.59.\h5.60.\h5.61.\h5.62.\h5.63.\h5.64.\h5.65.\h5.66.\h5.67.
\h5.68.\h5.69.\h5.70.\h5.71.\h5.72.\h5.73.\h5.74.\h5.75.\h5.76.\h5.77.\h5.78.
\h5.79.\h5.80.\h5.81.\h5.82.\h5.83.\h5.84.\h5.85.\h5.86.\h5.87.\h5.88.\h5.89.
\h5.90.\h5.91.\h5.92.\h5.93.\h5.94.\h5.95.\h5.97.\h5.98.\h5.99.\h5.101.
\h5.103.\h5.104.\h5.105.\h5.107.\h5.108.\h5.109.\h5.110.\h5.111.\h5.113.
\h5.115.\h5.116.\h5.117.\h5.119.\h5.121.\h5.122.\h5.123.\h5.125.\h5.127.
\h5.128.\h5.129.\h5.131.\h5.133.\h5.135.\h5.137.\h5.139.\h5.141.\h5.143.
\h5.145.\h5.149.\h5.153.\h5.161.\h5.165.\h6.26.\h6.27.\h6.28.\h6.30.\h6.32.
\h6.33.\h6.34.\h6.36.\h6.37.\h6.38.\h6.39.\h6.40.\h6.41.\h6.42.\h6.43.\h6.44.
\h6.45.\h6.46.\h6.47.\h6.48.\h6.49.\h6.50.\h6.51.\h6.52.\h6.53.\h6.54.\h6.55.
\h6.56.\h6.57.\h6.58.\h6.59.\h6.60.\h6.61.\h6.62.\h6.63.\h6.64.\h6.65.\h6.66.
\h6.67.\h6.68.\h6.69.\h6.70.\h6.71.\h6.72.\h6.73.\h6.74.\h6.75.\h6.76.\h6.77.
\h6.78.\h6.79.\h6.80.\h6.81.\h6.82.\h6.83.\h6.84.\h6.85.\h6.86.\h6.87.\h6.88.
\h6.89.\h6.90.\h6.91.\h6.92.\h6.93.\h6.94.\h6.95.\h6.96.\h6.97.\h6.98.\h6.99.
\h6.100.\h6.101.\h6.102.\h6.103.\h6.104.\h6.105.\h6.106.\h6.107.\h6.108.
\h6.109.\h6.110.\h6.111.\h6.112.\h6.114.\h6.115.\h6.116.\h6.117.\h6.118.
\h6.120.\h6.122.\h6.124.\h6.126.\h6.127.\h6.128.\h6.132.\h6.136.\h6.138.
\h6.140.\h6.142.\h6.144.\h6.148.\h6.150.\h6.152.\h6.153.\h6.156.\h6.160.
\h6.162.\h6.164.\h6.168.\h7.19.\h7.23.\h7.25.\h7.27.\h7.28.\h7.29.\h7.30.
\h7.31.\h7.32.\h7.33.\h7.34.\h7.35.\h7.36.\h7.37.\h7.38.\h7.39.\h7.40.\h7.41.
\h7.42.\h7.43.\h7.44.\h7.45.\h7.46.\h7.47.\h7.48.\h7.49.\h7.50.\h7.51.\h7.52.
\h7.53.\h7.54.\h7.55.\h7.56.\h7.57.\h7.58.\h7.59.\h7.60.\h7.61.\h7.62.\h7.63.
\h7.64.\h7.65.\h7.66.\h7.67.\h7.68.\h7.69.\h7.70.\h7.71.\h7.72.\h7.73.\h7.74.
\h7.75.\h7.76.\h7.77.\h7.78.\h7.79.\h7.80.\h7.81.\h7.82.\h7.83.\h7.84.\h7.85.
\h7.86.\h7.87.\h7.88.\h7.89.\h7.90.\h7.91.\h7.92.\h7.93.\h7.94.\h7.95.\h7.96.
\h7.97.\h7.98.\h7.99.\h7.100.\h7.101.\h7.102.\h7.103.\h7.104.\h7.105.\h7.106.
\h7.107.\h7.108.\h7.109.\h7.110.\h7.111.\h7.112.\h7.113.\h7.114.\h7.115.
\h7.116.\h7.117.\h7.119.\h7.120.\h7.121.\h7.122.\h7.123.\h7.124.\h7.125.
\h7.127.\h7.130.\h7.131.\h7.133.\h7.135.\h7.137.\h7.139.\h7.140.\h7.141.
\h7.142.\h7.143.\h7.145.\h7.147.\h7.149.\h7.151.\h7.153.\h7.154.\h7.155.
\h7.157.\h7.159.\h7.160.\h7.161.\h7.163.\h7.167.\h7.169.\h8.22.\h8.24.\h8.25.
\h8.26.\h8.27.\h8.28.\h8.29.\h8.30.\h8.31.\h8.32.\h8.33.\h8.34.\h8.35.\h8.36.
\h8.37.\h8.38.\h8.39.\h8.40.\h8.41.\h8.42.\h8.43.\h8.44.\h8.45.\h8.46.\h8.47.
\h8.48.\h8.49.\h8.50.\h8.51.\h8.52.\h8.53.\h8.54.\h8.55.\h8.56.\h8.57.\h8.58.
\h8.59.\h8.60.\h8.61.\h8.62.\h8.63.\h8.64.\h8.65.\h8.66.\h8.67.\h8.68.\h8.69.
\h8.70.\h8.71.\h8.72.\h8.73.\h8.74.\h8.75.\h8.76.\h8.77.\h8.78.\h8.79.\h8.80.
\h8.81.\h8.82.\h8.83.\h8.84.\h8.85.\h8.86.\h8.87.\h8.88.\h8.89.\h8.90.\h8.91.
\h8.92.\h8.93.\h8.94.\h8.95.\h8.96.\h8.97.\h8.98.\h8.99.\h8.100.\h8.101.
\h8.102.\h8.103.\h8.104.\h8.105.\h8.106.\h8.107.\h8.108.\h8.109.\h8.110.
\h8.111.\h8.112.\h8.113.\h8.114.\h8.115.\h8.116.\h8.118.\h8.119.\h8.120.
\h8.121.\h8.122.\h8.123.\h8.124.\h8.125.\h8.126.\h8.128.\h8.130.\h8.131.
\h8.132.\h8.134.\h8.136.\h8.137.\h8.138.\h8.140.\h8.142.\h8.143.\h8.144.
\h8.146.\h8.148.\h8.149.\h8.150.\h8.151.\h8.152.\h8.154.\h8.155.\h8.156.
\h8.158.\h8.160.\h8.161.\h8.162.\h8.163.\h8.164.\h8.166.\h8.167.\h8.168.
\h8.170.\h9.19.\h9.21.\h9.22.\h9.23.\h9.24.\h9.25.\h9.26.\h9.27.\h9.28.
\h9.29.\h9.30.\h9.31.\h9.32.\h9.33.\h9.34.\h9.35.\h9.36.\h9.37.\h9.38.\h9.39.
\h9.40.\h9.41.\h9.42.\h9.43.\h9.44.\h9.45.\h9.46.\h9.47.\h9.48.\h9.49.\h9.50.
\h9.51.\h9.52.\h9.53.\h9.54.\h9.55.\h9.56.\h9.57.\h9.58.\h9.59.\h9.60.\h9.61.
\h9.62.\h9.63.\h9.64.\h9.65.\h9.66.\h9.67.\h9.68.\h9.69.\h9.70.\h9.71.\h9.72.
\h9.73.\h9.74.\h9.75.\h9.76.\h9.77.\h9.78.\h9.79.\h9.80.\h9.81.\h9.82.\h9.83.
\h9.84.\h9.85.\h9.86.\h9.87.\h9.88.\h9.89.\h9.90.\h9.91.\h9.92.\h9.93.\h9.94.
\h9.95.\h9.96.\h9.97.\h9.98.\h9.99.\h9.100.\h9.101.\h9.102.\h9.103.\h9.104.
\h9.105.\h9.106.\h9.107.\h9.108.\h9.109.\h9.110.\h9.111.\h9.112.\h9.113.
\h9.114.\h9.115.\h9.116.\h9.117.\h9.118.\h9.119.\h9.120.\h9.121.\h9.122.
\h9.123.\h9.124.\h9.125.\h9.126.\h9.127.\h9.128.\h9.129.\h9.130.\h9.131.
\h9.132.\h9.133.\h9.134.\h9.135.\h9.136.\h9.137.\h9.138.\h9.139.\h9.141.
\h9.143.\h9.144.\h9.145.\h9.146.\h9.147.\h9.148.\h9.149.\h9.150.\h9.151.
\h9.153.\h9.154.\h9.155.\h9.156.\h9.157.\h9.158.\h9.159.\h9.161.\h9.162.
\h9.163.\h9.165.\h9.167.\h9.168.\h9.169.\h10.18.\h10.20.\h10.22.\h10.23.
\h10.24.\h10.25.\h10.26.\h10.27.\h10.28.\h10.29.\h10.30.\h10.31.\h10.32.
\h10.33.\h10.34.\h10.35.\h10.36.\h10.37.\h10.38.\h10.39.\h10.40.\h10.41.
\h10.42.\h10.43.\h10.44.\h10.45.\h10.46.\h10.47.\h10.48.\h10.49.\h10.50.
\h10.51.\h10.52.\h10.53.\h10.54.\h10.55.\h10.56.\h10.57.\h10.58.\h10.59.
\h10.60.\h10.61.\h10.62.\h10.63.\h10.64.\h10.65.\h10.66.\h10.67.\h10.68.
\h10.69.\h10.70.\h10.71.\h10.72.\h10.73.\h10.74.\h10.75.\h10.76.\h10.77.
\h10.78.\h10.79.\h10.80.\h10.81.\h10.82.\h10.83.\h10.84.\h10.85.\h10.86.
\h10.87.\h10.88.\h10.89.\h10.90.\h10.91.\h10.92.\h10.93.\h10.94.\h10.95.
\h10.96.\h10.97.\h10.98.\h10.99.\h10.100.\h10.101.\h10.102.\h10.103.\h10.104.
\h10.105.\h10.106.\h10.107.\h10.108.\h10.109.\h10.110.\h10.111.\h10.112.
\h10.113.\h10.114.\h10.115.\h10.116.\h10.117.\h10.118.\h10.119.\h10.120.
\h10.121.\h10.122.\h10.123.\h10.124.\h10.125.\h10.126.\h10.127.\h10.128.
\h10.129.\h10.130.\h10.131.\h10.132.\h10.133.\h10.134.\h10.135.\h10.136.
\h10.137.\h10.138.\h10.139.\h10.140.\h10.141.\h10.142.\h10.143.\h10.144.
\h10.145.\h10.146.\h10.148.\h10.149.\h10.150.\h10.151.\h10.152.\h10.153.
\h10.154.\h10.155.\h10.156.\h10.157.\h10.158.\h10.159.\h10.160.\h10.162.
\h10.163.\h10.164.\h10.166.\h10.167.\h10.168.\h10.169.\h10.170.
\end{picture}}

The advantage of this strategy is that it is easy to get a rather complete 
list of reflexive polytopes with a small number of points, which at the same 
time have a high probability for the existence of nef partitions and whose 
Hodge data are outside the range that is already completely covered by 
hypersurfaces.
Moreover, this class dissipates less time in computing the
nef partitions because of their small number of vertices.
Pursuing this strategy, a further step will be to increase the codimension by
one and to construct complete intersections using six-dimensional starting
polytopes with a small number of points.
\bigskip

{\it Acknowledgements.} This work was supported in part by the Austrian
        Research Funds FWF under grant Nr. P14639-TPH.

\newpage

\newpage

\begin{appendix}

\section{Hodge data}

\hspace{-.5cm}
\begin{tabular}{|r|r||r|r||r|r|r||r|r||r|r|}
\hline
$\mathrm{d_1}$ & $ \mathrm{w_{1_1},\dots,w_{n_1}} $ & $\mathrm{d_2}$ & $
\mathrm{w_{1_2},\dots,w_{n_2}} $ & $h^{11}$ & $h^{21}$ & $-\chi$ &
$\#\Delta\cap M$ & $\#\Delta^v $ & $\#\Delta^\ast\cap N $ &
$\#{\Delta^\ast}^v$\\\hline\hline
 3 & 1 1 1 0 0 0 0 &  4 & 0 0 0 1 1 1 1  & 2 &  59 & 114 &  350  & 12 &  8 & 7
\\\hline
 3 & 1 1 1 0 0 0 0 &  5 & 1 0 0 1 1 1 1  & 2 &  60 & 116 &  379  & 12 &  8 & 7
\\\hline
 3 & 1 1 1 0 0 0 0 &  4 & 0 0 0 1 1 1 1  & 2 &  62 & 120 &  350  & 12 &  8 & 7
\\\hline
3  & 1 1 1 0 0 0 0 &  6 & 0 1 1 1 1 1 1  & 2 &  62 & 120 &  381  & 12 &  8 & 7
\\\hline
 3 & 1 1 1 0 0 0 0 &  5 & 1 0 0 1 1 1 1  & 2 &  70 & 136 &  379  & 12 &  8 & 7
\\\hline
 3 & 1 1 1 0 0 0 0 &  6 & 0 1 1 1 1 1 1  & 2 &  76 & 148 &  381  & 12 &  8 & 7
\\\hline
 3 & 1 1 1 0 0 0 0 &  4 & 0 0 0 1 1 1 1  & 2 &  77 & 150 &  350  & 12 &  8 & 7
\\\hline
 3 & 1 1 1 0 0 0 0 &  8 & 3 0 0 1 1 2 1  & 2 & 100 & 196 &  496  & 16 &  9 & 8
\\\hline
 4 & 1 1 2 0 0 0 0 &  5 & 1 0 0 1 1 1 1  & 3 &  55 & 104 &  292  & 12 &  9 & 7
\\\hline
 4 & 1 1 2 0 0 0 0 &  6 & 1 0 0 1 1 1 2  & 3 &  55 & 104 &  282  & 12 &  9 & 7
\\\hline
 4 & 2 1 1 0 0 0 0 &  8 & 0 1 2 1 1 1 2  & 3 &  55 & 104 &  247  &  9 &  9 & 7
\\\hline
 4 & 2 1 1 0 0 0 0 & 12 & 0 1 3 2 2 2 2  & 3 &  55 & 104 &  265  &  9 &  9 & 7
\\\hline
 4 & 1 1 2 0 0 0 0 &  4 & 0 0 0 1 1 1 1  & 3 &  55 & 104 &  315  & 12 &  9 & 7
\\\hline
 3 & 1 1 1 0 0 0 0 &  5 & 0 0 0 1 1 1 2  & 3 &  56 & 106 &  340  & 18 &  9 & 8
\\\hline
 4 & 1 1 2 0 0 0 0 & 16 & 1 0 0 2 3 4 6  &13 &  15 &   4 &   117 & 20 & 15 &10
\\\hline
\end{tabular}

\vspace{0.5cm}

\noindent
Table 3. New Hodge numbers of toric CICYs using combined weight systems.\\

\noindent
\begin{tabular}{|r|r||r|r|r|r||r|r||r|r|}
\hline
$\mathrm{d}$ & $ \mathrm{w_{1},\dots,w_{n}} $ & $h^{11}$ & $h^{21}$ & $-\chi$ &
$\#\Delta\cap M$ & $\#\Delta^v $ & $\#\Delta^\ast\cap N $ &
$\#{\Delta^\ast}^v$\\\hline\hline
12 & 1 1 2 2 3 3    & 1 &  61 & 120 & 407 &  6 &  7 & 6\\\hline
 6 & 1 1 1 1 1 1    & 1 &  73 & 144 & 462 &  6 &  7 & 6\\\hline
 8 & 1 1 1 1 2 2    & 1 &  73 & 144 & 483 &  6 &  7 & 6\\\hline
 6 & 1 1 1 1 1 1    & 1 &  89 & 176 & 462 &  6 &  7 & 6\\\hline
12 & 1 2 2 2 2 3    & 2 &  62 & 120 & 321 &  6 &  8 & 6\\\hline
 9 & 1 1 1 2 2 2    & 2 &  68 & 132 & 434 & 12 &  8 & 7\\\hline
10 & 1 1 2 2 2 2    & 2 &  68 & 132 & 378 &  6 &  8 & 6\\\hline
\end{tabular}

\vspace{0.5cm}

\noindent
Table 4. New Hodge numbers of toric CICYs using one weight system.\\

\noindent
\begin{tabular}{|r|r||r|r||r|r|r||r|r||r|r|}
\hline
$\mathrm{d}$ & $ \mathrm{w_{1},\dots,w_{n}} $ & $\mathrm{d_1}$ & $\mathrm{d_2}$
&  $h^{11}$ & $h^{21}$ & $-\chi$ & $\#\Delta\cap M$ & $\#\Delta^v $ &
$\#\Delta^\ast\cap N $ & $\#{\Delta^\ast}^v$\\\hline\hline
12 & 1 1 2 2 3 3 & 6 & 6 & 1 &  61 & 120 & 407 & 6  & 7  & 6\\\hline
 6 & 1 1 1 1 1 1 & 2 & 4 & 1 &  73 & 144 & 462 & 6  & 7  & 6\\\hline
 8 & 1 1 1 1 2 2 & 4 & 4 & 1 &  73 & 144 & 483 & 6  & 7  & 6\\\hline
 6 & 1 1 1 1 1 1 & 2 & 4 & 1 &  89 & 176 & 462 & 6  & 7  & 6\\\hline
 8 & 1 1 1 1 1 3 & 4 & 4 & 1 & 129 & 256 & 636 & 10 & 8  & 7\\\hline
12 & 1 2 2 2 2 3 & 6 & 6 & 2 &  62 & 120 & 321 & 6  & 8  & 6\\\hline
 9 & 1 1 1 2 2 2 & 4 & 5 & 2 &  68 & 132 & 434 & 12 & 8  & 7\\\hline
10 & 1 1 2 2 2 2 & 4 & 6 & 2 &  68 & 132 & 378 & 6  & 8  & 6\\\hline
14 & 1 2 2 3 3 3 & 8 & 6 & 3 &  47 &  88 & 294 &12  & 9  & 7\\\hline
16 & 1 2 2 3 4 4 & 8 & 8 & 3 &  55 & 104 & 327 & 8  & 9  & 7\\\hline
\end{tabular}

\vspace{0.5cm}

\noindent
Table 5. New Hodge numbers of toric CICYs using Minkowski sums.\\
\end{appendix}

\bye